\documentclass[12pt]{article}
\usepackage{amssymb}
\usepackage{epsfig}

\newcommand{\triangleboard}{\textsf{Triangle}}
\newcommand{\rhombusboard}{\textsf{Rhombus}}

\newenvironment{packed_enumerate}{
\setlength{\parsep}{0pt}
\setlength{\parskip}{0pt}
\begin{enumerate}
  \setlength{\itemsep}{1pt}
  \setlength{\parsep}{0pt}
  \setlength{\parskip}{0pt}
}{\end{enumerate}}

\begin{document} 

\begin{center}
\vspace*{0.1in}
{\huge Diamond Solitaire}
\vskip 20pt
{\bf George I. Bell}\\
September 2005\footnote{Original version at {\tt http://gpj.connectfree.co.uk/gpjw.htm} \\
\hspace*{0.25in}Converted to \LaTeX ~by the author with some modifications to the text, November 2007.} \\
{\tt gibell@comcast.net}\\
\end{center}
\vskip 30pt 
\centerline{\bf Abstract}
\noindent
We investigate the game of peg solitaire on different board shapes,
and find those of diamond or rhombus shape have interesting properties.
When one peg captures many pegs consecutively,
this is called a sweep.
Rhombus boards of side 6 have the property that no matter which peg
is missing at the start, the game can be solved to one peg
using a maximal sweep of length 16.
We show how to construct a solution on a rhombus board of side $6i$,
where the final move is a maximal sweep of length $r$,
where $r=(9i-1)(3i-1)$ is a ``rhombic matchstick number".

\pagestyle{myheadings}
\markright{The Games and Puzzles Journal---Issue 41, September-October 2005\hfill}

\thispagestyle{empty} 
\baselineskip=15pt 
\vskip 30pt 

\section{Introduction}

\begin{figure}[htb]
\centering
\epsfig{file=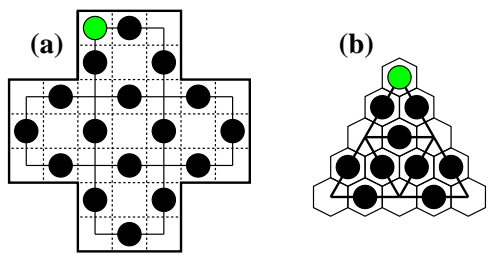}
\caption{Maximal sweeps on the standard boards, only the second is a super-sweep.
The special green (or shaded) peg performs the sweep.
(a) A 16-sweep on the 33-hole cross-shaped board.
(b) A 9-sweep on the 15-hole triangular board, \triangleboard$(5)$.}
\label{fig1}
\end{figure}

\noindent
Peg solitaire is a one-person game usually played on a
33-hole cross-shaped board,
or a 15-hole triangular board (Figure~\ref{fig1}).
In the first case the pattern of holes (or board locations)
come from a square lattice,
while in the second case they come from a triangular lattice.
The usual game begins with the board filled by pegs
(or marbles) except for a single board location, left vacant.
The player then jumps one peg over another into
an empty hole, removing the jumped peg from the board.
The goal is to choose a series of jumps that finish with one peg.
The general problem of going from a board position with one peg missing
to one peg will be called a \textbf{peg solitaire problem}.
If the missing peg and final peg are in the same place,
we call it a \textbf{complement problem}, because the starting
and ending board positions are complements of one another
(where every peg is replaced by an empty hole and vice versa).

\noindent
Initially, we will consider peg solitaire on boards of
rather arbitrary shape.
The basic problem may seem hard enough, but we will add
more conditions or constraints on the solution.
The reason for this may not be immediately clear,
but adding constraints to a problem with many solutions
can make it easier to find one, and the solutions themselves
can be quite remarkable.
The approach also brings out certain board shapes with
interesting properties.

\noindent
Suppose we add the constraint that the peg solitaire problem must
finish in the most dramatic way possible: with one peg ``sweeping off"
all the remaining pegs, and finishing as the sole survivor.
Here we need to distinguish between a \textbf{jump}---one
peg jumping over another, and a \textbf{move}---one or more
jumps by the same peg.
Any single move that captures $n$ pegs is
a \textbf{sweep} of length $n$, or an \textbf{\textit{n}-sweep}.
In this terminology we want to finish a peg solitaire problem
with the longest sweep possible.

\noindent
A sweep that has the longest length geometrically possible
on a board is called a \textbf{maximal} sweep.
For most board shapes, maximal sweeps cannot be reached
in the solution to a peg solitaire problem.
In GPJ \#36 \cite{GPJ36}, we considered the triangular board
of side $n$, called \triangleboard$(n)$.
We saw that for $n$ odd, the board \triangleboard$(n)$
supports an even more remarkable sweep that is maximal and also has
the property that every hole in the board is affected by the sweep
in that it is either jumped over or is the starting or ending hole of some
jump in the sweep.
Such remarkable sweeps need a special name---we call them
\textbf{super-sweeps}.

\noindent
In Figure~\ref{fig1}a, the central hole is not touched or
jumped over by the sweep, thus while this sweep is the longest possible
it is not a super-sweep, and a super-sweep is not possible on this board.
It's not hard to see that super-sweeps are never possible on square
lattice boards, except for certain trivial or degenerate cases
(like a 1-dimensional board).
Non-trivial super-sweeps are only possible on a triangular lattice, an
example is Figure~\ref{fig1}b.

\noindent
Neither of the sweep patterns of Figure~\ref{fig1} can be reached during
the solution to a peg solitaire problem on that board.
The easiest way to see this is to try to figure out how one
could have arrived at the sweep board position, or equivalently
try to play \textit{backwards} from the sweep position.
In GPJ \#36 \cite{GPJ36},
we saw in the \textbf{Forward/Backward Theorem}
that backward play is equivalent to forward play from the
complement of the board position.
This concept is referred to as the ``time reversal trick" in
Winning Ways For Your Mathematical Plays \cite{WinningWays}.

\noindent
If we take the complement of either board position in
Figure~\ref{fig1}, we find that no jump is possible.
This proves that there is no ``parent board position" from which
the sweep position can arise, it cannot occur during the solution
to a peg solitaire problem.
If we take the complement of any super-sweep pattern,
and a jump is possible, this means there would have to be
two consecutive empty holes
in the original super-sweep.
Because it is a super-sweep, \textit{both} of these holes must be
the starting or ending locations of some jump in the sweep,
which is impossible.
Therefore a super-sweep \textit{can never occur} in the
solution to a peg solitaire problem.

\noindent
If super-sweeps cannot occur in peg solitaire problems, the reader
may wonder why we are wasting our time with them.
The answer is provided by the triangular boards and GPS \#36 \cite{GPJ36}:
the super-sweep pattern of Figure~\ref{fig1}b \textit{can} be
reached in a problem on \triangleboard$(6)$.
This 9-sweep is no longer a super-sweep
with respect to \triangleboard$(6)$,
\textit{but it is still a maximal sweep}.
Loosely speaking, a super-sweep may still be reachable in a
peg solitaire problem on a board one size larger\footnote{Amazingly,
the 16-sweep of Figure~\ref{fig1}a can also be reached
in a peg solitaire problem on a cross-shaped board one size larger,
the 45-hole ``Wiegleb's board" \cite{BGS}.}.

\section{General boards with super-sweeps}

\begin{figure}[htb]
\centering
\epsfig{file=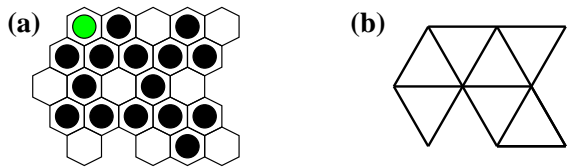}
\caption{An 8-sided (non-convex) polygon board and the associated super-sweep graph.
(a) An unusual 24-hole board, set up in the board position of a hypothetical super-sweep
(of length 15).
(b) The graph of the hypothetical super-sweep.}
\label{fig2}
\end{figure}

\noindent
Let us consider peg solitaire boards on a triangular lattice
where the board shape is a general polygon.
On a triangular lattice the corners of the board
are restricted to multiples of $60^\circ$.
For which such boards is it possible to have a super-sweep?
It is clear that the board edges must all have an odd length (in holes),
because the super-sweep must pass through all the corners of the board
(by definition, a \textbf{corner} hole cannot be jumped over,
and consequently the super-sweep must jump into it).
Consider, for example, the 8-sided polygon board of Figure~\ref{fig2}.

\noindent
Can the green peg (top left) in Figure~\ref{fig2}a make a
tour of the board, sweeping off all 15 remaining pegs?
We can answer this question by looking at the graph formed by
the sweep, shown separately in Figure~\ref{fig2}b
rather than on top of the board as in Figure~\ref{fig1}.
In the language of graph theory the super-sweep is an Euler path
on this graph, i.e. a path that traverses every edge exactly once.
One of the most basic theorems of graph theory states that
an Euler path is possible if and only if there
are either zero or two nodes of odd degree.
Looking at Figure~\ref{fig2}b,
we see that there are four nodes of odd degree,
hence an Euler path is not possible.
The board in Figure~\ref{fig2} \textit{does not} support a super-sweep.

\noindent
Using the Euler path theorem we can eliminate many boards
from having super-sweeps.
Hexagonal or star-shaped boards have hexagonal symmetry,
but do not support super-sweeps because the associated
sweep graphs have six nodes of odd degree.
If we restrict ourselves to board shapes that are \textit{convex polygons}
(no $240^\circ$ or $300^\circ$ corners) things are particularly simple.
For a convex board on a triangular lattice, only
three board shapes can have a super-sweep:
\begin{packed_enumerate}
\item \textbf{Triangles}, and only equilateral triangles
are possible on a triangular lattice.
Because every node in the super-sweep graph has even degree,
the super-sweep always begins and ends at the
same board location.
\item \textbf{Parallelograms}, which have alternating $60^\circ$ and
$120^\circ$ corners as you go around the circumference.
The super-sweep must begin at one $120^\circ$ corner and
end at the other.
\item \textbf{Trapezoids}, which have two $60^\circ$ corners followed
by two $120^\circ$ corners.
The super-sweep must begin at one $120^\circ$ corner
and end at the other.
These boards can also be considered as
equilateral triangles with one corner cut off.
\end{packed_enumerate}

\begin{figure}[htb]
\centering
\epsfig{file=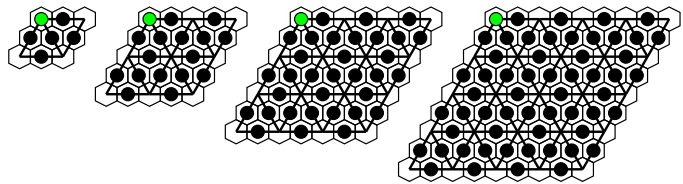}
\caption{Super-sweeps on \rhombusboard$(n)$ where $n=3$, $5$, $7$ and $9$.
These sweeps have lengths of 5, 16, 33, and 56 respectively,
and end at the lower right corner.}
\label{fig3}
\end{figure}

\noindent
A rhombus is a special case of a parallelogram where all
four sides have the same length $n$,
we'll call this board \rhombusboard$(n)$.
By rotating them $60^\circ$, they become diamond-shaped
and could also be called \textsf{Diamond} boards.
In the remainder of this paper, we'll go into some of the remarkable
properties of these boards.
Figure~\ref{fig3} shows the first four rhombus board super-sweeps.

\noindent
The length of this sweep, for odd $n$, is $(3n+1)(n-1)/4$.
Since the total number of holes in the board is $n^2$,
this sweep removes nearly $3/4$ of the pegs on the board.
For $n$ = $3$, $5$, $7$, $9$, $11$, $13$, $\ldots$,
the sequence of sweep lengths is $5$, $16$, $33$, $56$, $85$, $120$, $\ldots$,
called the ``rhombic matchstick sequence" \cite{OEIS}
because it is the number of matchsticks needed to construct a rhombus
(with $(n-1)/2$ matchsticks on a side).

\section{\rhombusboard$(6)$}

\noindent
This 36-hole board has several unusual properties.
It is also of a reasonable size for playing by hand,
and for computational searches.
This board is equivalent to the 6x6 square board on a square lattice,
\textit{with the addition of moves along one diagonal}.
It is therefore possible to play this board using a chess or go board,
although this is not recommended because the symmetry of the board
is obscured.
For playing by hand I recommend using part of a Chinese Checkers board.
The ideal board for playing by hand is a computer \cite{Bellweb},
because we can easily take back moves and record move sequences.

\noindent
The board \rhombusboard$(6)$ is a \textbf{null-class} board.
For a definition of this term, see \cite{Beasley} or \cite{Bellweb},
the important concept is that only on null-class boards can
a complement problem be solvable.
\rhombusboard$(6)$ is the smallest rhombus board on which a complement
problem is solvable, and in fact \textit{all} complement problems
are easily solvable.

\begin{figure}[htb]
\centering
\epsfig{file=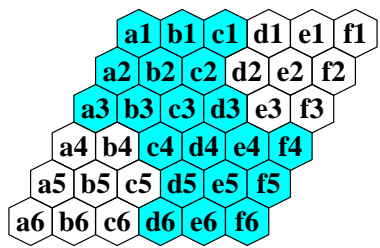}
\caption{\rhombusboard$(6)$ hole coordinates.
Potential finishing locations for a problem including a maximal sweep
(16-sweep) are shaded blue.}
\label{fig4}
\end{figure}

\noindent
The longest sweep geometrically possible on \rhombusboard$(6)$
has length 16 (as in Figure~\ref{fig3}).
Can 16-sweeps occur in solutions to peg solitaire problem
on this board?
Here we aren't limiting the 16-sweep to be the last move,
but leave open the possibility that it could happen at any move.
We note that there are only a few places where the 16-sweep
can begin and end.
It can go from a1 to e5, b2 to f6, b1 to f5, or any symmetric variant of these.
The 16-sweep can be the final move, or it can be the second
to the last move, for example the 16-sweep can go from a1 to e5,
followed by f6-d4, or e6-e4\footnote{Notation: moves are denoted by the starting
and ending coordinates of the jumps, separated by dashes.}.
The 16-sweep can even be the 3rd to the last move,
from a1 to e5, followed by e6-e4 and f5-d3.

\noindent
In Figure~\ref{fig4}, all potential finishing locations
for solutions containing a 16-sweep are shaded blue.
Not all are feasible,
the finishing 16-sweep from b1 to f5 is in fact impossible to reach
from any starting vacancy, as discovered by computational search.
However all other configurations of the 16-sweep can be reached.
In fact, starting from any vacancy on the board, there is a solution
with a maximal sweep (16-sweep) that finishes with one peg.
This board is the only one we know of, on a square or triangular lattice,
with this amazing property.

\noindent
Another unusual property of \rhombusboard$(6)$ is that solutions to
complement problems can include maximal sweeps.
Here are four different complement problems that
can be solved using a maximal sweep:
\begin{packed_enumerate}
\item e5 complement: solve with the last move a 16-sweep.
\item d4 complement: solve with the second to last move a 16-sweep.
\item e4 complement: solve with the second to last move a 16-sweep.
\item d3 complement: solve with the third to last move a 16-sweep.
\end{packed_enumerate}

\begin{figure}[htb]
\centering
\epsfig{file=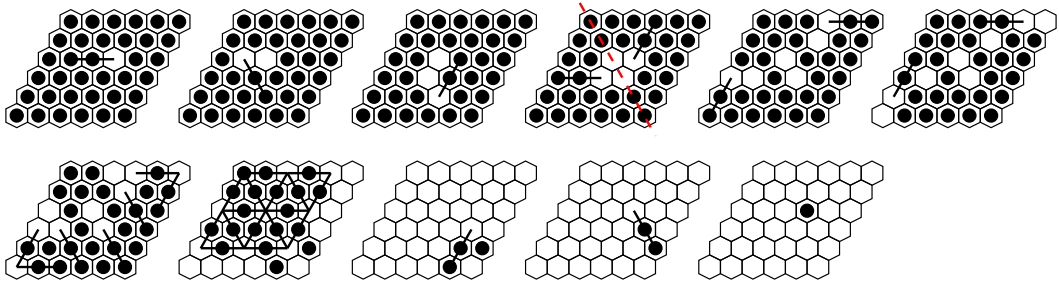}
\caption{A symmetric solution to the d3 complement (problem \#4).
This solution has 17 moves.
Note that more than one move is sometimes shown between board snapshots.}
\label{fig5}
\end{figure}

\noindent
These problems are most easily solved by attempting to play backward,
or equivalently by playing forward from the complement of the
board position before the 16-sweep.
All four problems make good challenges to solve by hand;
they are easy to solve using a computer
(provided you don't try to solve them playing forward).
In Figure~\ref{fig5} we show a solution to problem \#4.
This solution is interesting in that after the third move, the board position
is symmetric about the red dashed line.
After that moves are done in pairs, or are themselves symmetric,
preserving the symmetry up until the last two moves.

\noindent
Any peg solitaire problem on this board begins with 35 pegs and
finishes with one, so a solution consists of exactly 34 jumps.
The number of moves, however, can be less than this, and
an interesting problem is to find solutions in as few moves
as possible.
This is different from finding solutions
with maximal sweeps, and answers are
more difficult to obtain.
The minimal solution length can be found
using computational search methods \cite{BellDiag,BellTriang}.

\begin{figure}[htb]
\centering
\epsfig{file=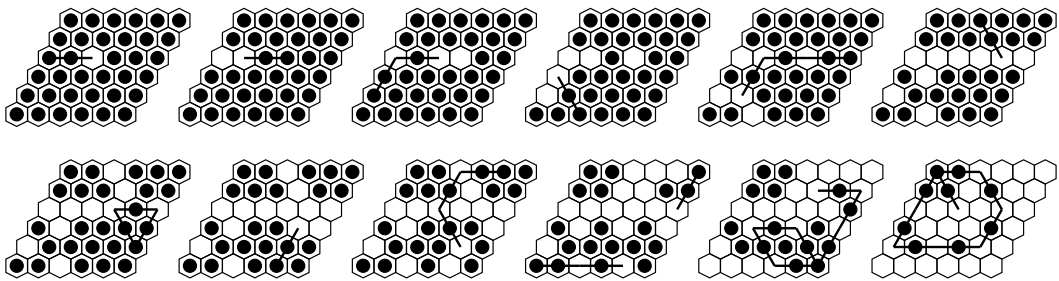}
\caption{A 13-move solution to the c3-complement.
The last four moves originate from corners, an unusual property for a 13-move solution.
Note that more than one move is sometimes shown between board snapshots.}
\label{fig6}
\end{figure}

\begin{figure}[htb]
\centering
\epsfig{file=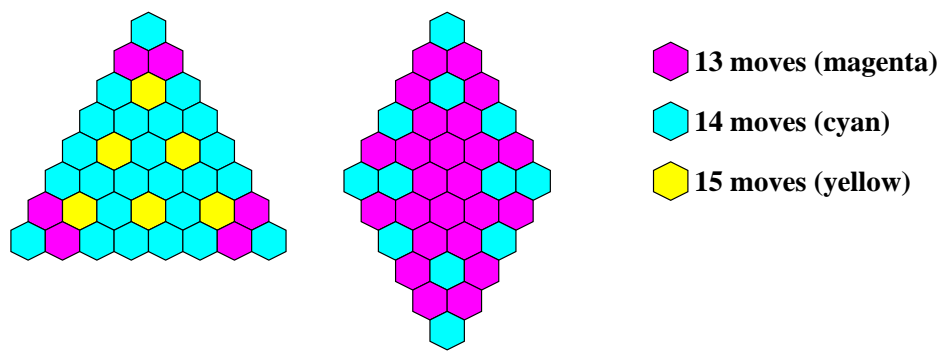}
\caption{Colors indicate the length of the shortest solution to a
complement problem on \triangleboard$(8)$ and \rhombusboard$(6)$
(magenta = 13 moves, cyan = 14 moves, yellow = 15 moves).
The \rhombusboard$(6)$ board has been rotated to its ``diamond"
configuration to show the symmetry.}
\label{fig7}
\end{figure}

\noindent
If we take into account all possible starting and finishing locations
for a peg solitaire problem on this board, we find there are 120
distinct problems.
I have only solved the complement problems, for
there are only 12 of them.
Of these 12, I have found that 7 can be solved
in a minimum of 13 moves, with the rest requiring 14 moves
(see Figure~\ref{fig7} for all results).
A sample 13-move solution is shown in Figure~\ref{fig6}.

\noindent
Using the ``Merson region" analysis of GPJ \#36 \cite{GPJ36},
it is possible to prove that the 13-move solution in Figure~\ref{fig6}
is the shortest possible.
In general, however, we rely on computational search
to establish the minimum.
For more information on minimal length solutions on
\rhombusboard$(6)$ and other boards see my web site \cite{Bellweb}.

\noindent
The smallest integer (greater than one) that is
a triangular number and a perfect square is 36.
This is reflected in the fact that \triangleboard$(8)$ and
\rhombusboard$(6)$ both have 36 holes.
Because of this, it is interesting to compare the
properties of these two boards.
Figure~\ref{fig7} shows the minimum length solution of a complement problem
by color for each of these boards.
Note that \rhombusboard$(6)$ in general supports slightly
shorter solutions, with none requiring 15 moves.

\section{Maximal sweeps on \rhombusboard$(6i)$}

\noindent
In GPJ \#36 \cite{GPJ36}, we found that maximal sweeps on
\triangleboard$(6)$ and \triangleboard$(8)$ could occur in
peg solitaire games on these boards.
Although a proof has not been found, computational results
suggest that maximal sweeps cannot be reached on any larger
triangular boards.

\noindent
We have just shown that a maximal sweep can be reached by a peg solitaire
problem on \rhombusboard$(6)$, but what about larger rhombus boards?
One might suspect that maximal sweeps would eventually become unreachable,
as with the triangular boards.
Somewhat remarkably, however, this is not the case.

\noindent
{\bf Theorem} For any $i>0$, there exists a solution
to a peg solitaire problem on \rhombusboard$(6i)$
where the last move is a maximal sweep of length $(9i-1)(3i-1)$.

\begin{figure}[htb]
\centering
\epsfig{file=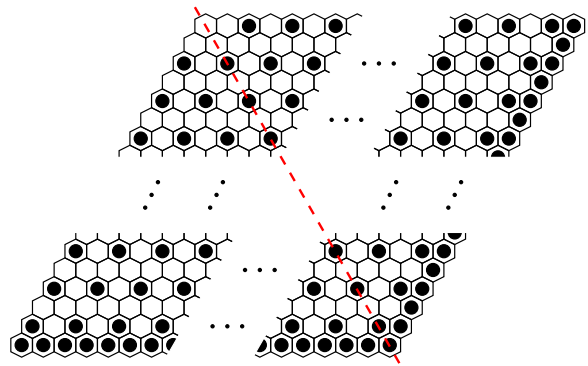}
\caption{The complement of the sweep pattern on \rhombusboard$(6i)$.
Note that the board position is symmetric about the red dashed line.}
\label{fig8}
\end{figure}

\noindent
{\it Proof:} By the \textbf{Forward/Backward Theorem} \cite{GPJ36},
it suffices to show that the complement of some
maximal sweep pattern can be reduced to one peg.
The sweep pattern we choose begins at the upper left corner, and
ends one hole up and left from the lower right corner.
When we take the complement,
this results in the board position of Figure~\ref{fig8}.

\noindent
The case $i=1$, or \rhombusboard$(6)$,
has already been solved (problem \#1 in the previous section).
We will use induction to prove the general case,
starting with $i=2$.
The solution proceeds through three phases ($A$, $B$ and $C$).
We apply the moves of phase $A$ once, then $B$ $(i-2)$ times,
followed by phase $C$ once.

\begin{figure}[htbp]
\centering
\epsfig{file=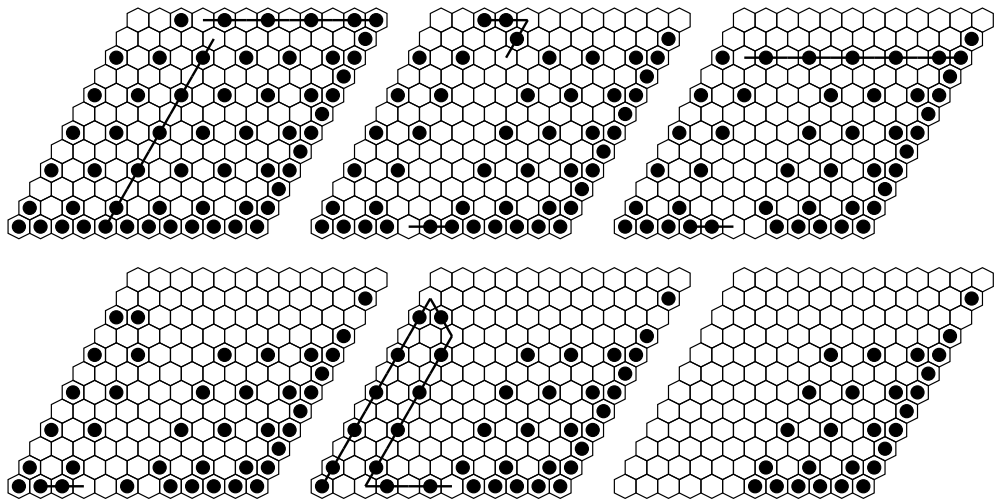}
\caption{The eight moves of phase $A$.}
\label{fig9}
\end{figure}

\noindent
Phase $A$, shown in Figure~\ref{fig9},
consists of eight moves that clear out the leftmost 6 columns
of the board and the upper 4 rows,
except for the last (rightmost) column.
If the board is larger than the \rhombusboard$(12)$ shown
in Figure~\ref{fig9}, the long multi-jump moves must be extended accordingly.
We are left with a very similar board pattern as the one
we started with, just reduced in size.
Note, however, that the final board position is no longer symmetric.

\begin{figure}[htbp]
\centering
\epsfig{file=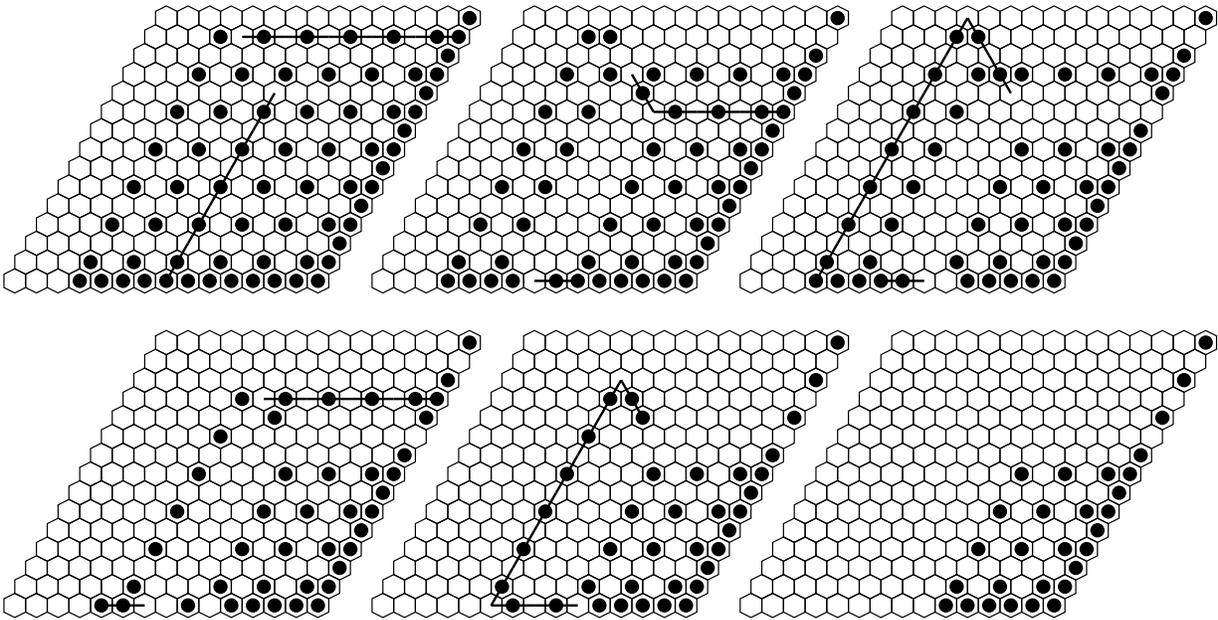}
\caption{The nine moves of phase $B$.
The moves here are shown on \rhombusboard$(15)$ to save space.
Phase $B$ will actually be applied to boards at least as large as
\rhombusboard$(18)$.}
\label{fig10}
\end{figure}

\noindent
Phase B, shown in Figure~\ref{fig10},
is nine moves that reduce the sweep pattern
by 6 rows and 6 columns.
As before if the board is larger than shown the multi-jump moves
are extended.
After applying phase $B$ $j$ times the leftmost $6j+6$ columns will
be empty, and the topmost $6j+4$ rows will also be empty, except for a
trail of pegs in the rightmost column that will be taken by the final move.

\begin{figure}[htbp]
\centering
\epsfig{file=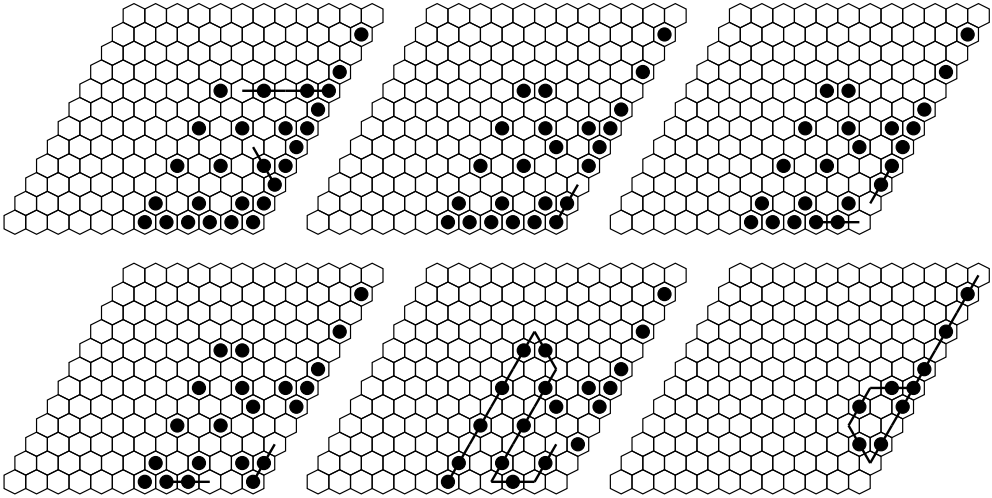}
\caption{The nine moves of phase $C$ (also possible in eight moves).}
\label{fig11}
\end{figure}

\noindent
Finally, phase $C$ is executed to take the board down to a one peg
in the upper right corner.
The nine moves of this phase are shown in Figure~\ref{fig11}.
Only in phase $C$ is the fact that the side is divisible by 6 needed.
For only on such boards can the final peg finish in the
upper right hand corner.

\noindent
Putting together all three phases,
it only takes $9i-1$
moves to clear the sweep pattern of Figure~\ref{fig8}.
To find the solution ending in the maximal sweep, we begin from a vacancy
at the upper right corner, and execute the jumps of
phases $A$, $B$ and $C$ in exactly the reverse order.
Thus we execute the jumps in phase $C$ reversed, followed by
$(i-2)$ phase $B$'s, and then phase $A$, all in reverse order.
The long sweeps become individual jumps,
and the solution ending in the maximal sweep
has significantly more than
$9i-1$ moves.

\begin{figure}[htbp]
\centering
\epsfig{file=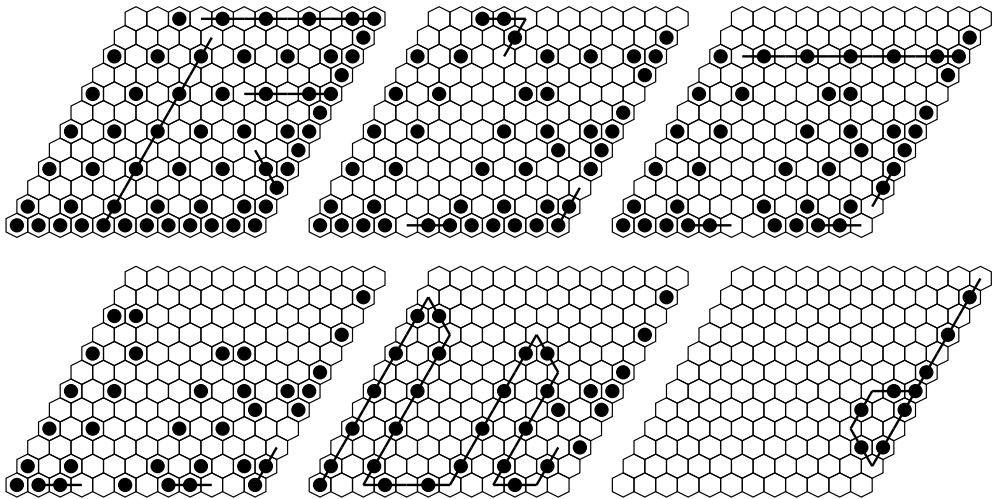}
\caption{Putting phases $A$ and $C$ together for a backward
solution on \rhombusboard$(12)$.}
\label{fig12}
\end{figure}

\noindent
In Figure~\ref{fig12} on \rhombusboard$(12)$,
we show the complete solution that reduces
the complement of the sweep pattern to one peg.
In this case we need only execute phases $A$ and $C$,
and in Figure~\ref{fig12} the
two phases have been interleaved and are no longer visually separate.
This reduces the number of diagrams to show the solution, but the
inductive nature of the solution becomes much harder to see.
It is unfortunate that a Chinese Checkers board is too small to play
this solution on.
If you can find a large enough board, it is interesting to play the
moves in this solution in exactly the reverse order,
and watch as the sweep position magically appears.
The final sweep in the reversed solution has length 85.

\noindent
An integer that is a perfect square and a triangular number
is called a \textbf{square triangular number} \cite{OEIS}.
As was the case with \rhombusboard$(6)$ and \triangleboard$(8)$,
which both have 36 holes, each square triangular number corresponds
to a rhombus and triangular board of the same size.
If the side of the rhombus board is divisible by 6, and the side of the
triangular board by 12, then both have long sweep finishes by the
above analysis and GPJ \#36 \cite{GPJ36}.

\noindent
After the 36-hole boards, the next time this occurs is with
\rhombusboard$(204)$ and \triangleboard$(288)$,
which both have $41,616$ holes.
By our inductive arguments, we can construct solutions to
peg solitaire problems on these boards that finish
with sweeps of length $30,805$ and $30,793$, respectively.
The next larger such boards are \rhombusboard$(235,416)$ and
\triangleboard$(332,928)$, boards with over 55 billion holes,
that can finish with sweeps over 41 billion in length.

\section{Conclusions}

\noindent
Peg solitaire on a triangular lattice is a fascinating game, and one
that has not been studied to the extent that ``normal"
(square lattice) peg solitaire has.
This paper, and GPJ \#36 \cite{GPJ36} have shown that triangular
lattice peg solitaire is well suited for inductive arguments.
Inductive arguments have also been used to create an algorithm
for triangular boards of any size that can reduce any (solvable)
single vacancy down to one peg \cite{BellTriang}.
I believe inductive arguments are possible for square lattice boards,
but the reduction from 6 jump directions to 4 seems to make such
arguments more difficult.

\noindent
While the results in Section~3 for \rhombusboard$(6)$
were obtained by computational search,
the long sweep finishes on \rhombusboard$(6i)$
were found by hand.
The computer was still of significant help, but only in providing an
interface to play the game on the large boards required.

\noindent
There remain many unanswered questions regarding rhombus and triangular boards.
Can maximal sweeps be reached on peg solitaire problems
on \rhombusboard$(2i)$?
I have been able to answer this question in the affirmative for $2\le i\le 9$.
It would also be nice to prove that maximal sweeps are not reachable
in peg solitaire problems on triangular boards larger than
\triangleboard$(8)$ (or find a counterexample).

\end{document}